\pdfoutput=1 

\documentclass[11pt]{article}

\usepackage{graphicx}
\usepackage{amssymb}
\usepackage{JASA_manu}

\usepackage{amsgen,amsmath,amstext,amsbsy,amsopn,amssymb}
\usepackage{epsfig}
\usepackage{cases}
\usepackage{amscd}
\usepackage{graphicx}% Include figure files
\usepackage{dcolumn}% Align table columns on decimal point
\usepackage{bm}% bold math
\usepackage{subfigure}
\usepackage{color}
\usepackage{enumerate}
\usepackage[round]{natbib}

\numberwithin{equation}{section}
\newtheorem{thm}{Theorem}[section]

\newtheorem{prop}{Proposition}[section]

\newcommand {\reals}  {\mathbb{R}}
\newcommand{\x}{{\mathbf{x}}}

\newcommand{\X}{{\mathbf{X}}}

\newcommand{\tx}{{\mathbf{\tilde{x}}}}

\newcommand{\bSigma}{\boldsymbol{\Sigma}}

\begin{document}

\title{Uniform Correlation Mixture of Bivariate Normal Distributions and Hypercubically-contoured Densities That Are Marginally Normal}
\author{Kai Zhang \\
Department of Statistics and Operations Research \\ University of North Carolina, Chapel Hill, 27599\\
email: \texttt{zhangk@email.unc.edu} \\
%\vskip 1in
Lawrence D. Brown \\
Department of Statistics, The Wharton School\\ University of Pennsylvania, Philadelphia, 19104\\
email: \texttt{lbrown@wharton.upenn.edu} \\
%\vskip .1in
Edward George \\
Department of Statistics, The Wharton School\\ University of Pennsylvania, Philadelphia, 19104\\
email: \texttt{edgeorge@wharton.upenn.edu} \\
%\vskip .1in
Linda Zhao \\
Department of Statistics, The Wharton School\\ University of Pennsylvania, Philadelphia, 19104\\
email: \texttt{lzhao@wharton.upenn.edu} \\
}

\maketitle

\newpage

\mbox{}
\vspace*{2in}
\begin{center}
\textbf{Author's Footnote:}
\end{center}
Kai Zhang is Assistant Professor, Department of Statistics and Operations Research, the University of North Carolina, Chapel Hill (E-mail: zhangk@email.unc.edu). Lawrence D. Brown, Edward George, and Linda Zhao are Professors, Department of Statistics, the Wharton School, the University of Pennsylvania. We appreciate insightful comments and suggestions from Robert Adler, Barry Arnold, James Berger, Richard Berk, Shankar Bhamidi, Peter Bickel, Amarjit Budhiraja, Andreas Buja, Dennis Cook, Bradley Efron, Jianqing Fan, Andrew Gelman, Ben Goodrich, Shane Jensen, Tiefeng Jiang, Iain Johnstone, Abba Krieger, J. Steve Marron, Henry McKean, Xiao-Li Meng, Dylan S. Small, J. Michael Steele, Yazhen Wang, and Cun-Hui Zhang. We also thank the editors and reviewers for careful reading and helpful comments. This work is supported by grants DMS-1309619 and DMS-1310795 from the U.S. National Science Foundation.

\newpage
\begin{center}
\textbf{Abstract}
\end{center}
The bivariate normal density with unit variance and correlation $\rho$ is well-known. We show that by integrating out $\rho$, the result is a function of the maximum norm. The Bayesian interpretation of this result is that if we put a uniform prior over $\rho$, then the marginal bivariate density depends only on the maximal magnitude of the variables. The square-shaped isodensity contour of this resulting marginal bivariate density can also be regarded as the equally-weighted mixture of bivariate normal distributions over all possible correlation coefficients. This density links to the Khintchine mixture method of generating random variables. We use this method to construct the higher dimensional generalizations of this distribution. We further show that for each dimension, there is a unique multivariate density that is a differentiable function of the maximum norm and is marginally normal, and the bivariate density from the integral over $\rho$ is its special case in two dimensions.

\vspace*{.3in}


\newpage

\section{Introduction}
It is well-known that a multivariate distribution that has normal marginal distributions is not necessarily jointly multivariate normal (in fact, not even when the distribution is conditionally normal, see \citet{gelman1991}), i.e., a $p$-dimensional multivariate distribution $\X=(X_1,\ldots,X_p)$ that has marginal standard normal densities $\phi(x_1), \phi(x_2), \ldots, \phi(x_p)$ and marginal distribution functions $\Phi(x_1), \Phi(x_2), \ldots, \Phi(x_p)$ may not have a jointly multivariate normal density
\begin{equation}
  f(x_1,\ldots,x_p|\bSigma) = (2\pi)^{-p/2}|\bSigma|^{-1/2} \exp\large\{-\x^T\bSigma^{-1}\x/2\large\}
\end{equation}
for some $p$ by $p$ correlation matrix $\bSigma.$ Classical examples of such distributions can be found in \citet{feller1971} and \citet{kotz2004}.

In this paper we focus on one particular class of such distributions that arises from uniform mixtures of bivariate normal densities over the correlation matrix. When $p=2,$ the bivariate normal density with unit variances is well-known:
\begin{equation}
  f(x_1,x_2|\rho) = {1 \over 2 \pi \sqrt{1-\rho^2}} \exp \left\{-{x_1^2+x_2^2-2\rho x_1x_2 \over 2(1-\rho^2)}\right\}
\end{equation}
for some $-1 \le \rho \le 1.$ By a uniform correlation mixture of bivariate normal density we mean the bivariate density function $f(x_1,x_2)$ below:
\begin{equation}
f(x_1,x_2)=\int_{-1}^1{1 \over 2}f(x_1,x_2|\rho) d \rho.
\end{equation}
This type of continuous mixture of bivariate normal distributions has been used in applications such as imaging analysis (\citet{Aylward1997}). We show that such a uniform correlation mixture results in a bivariate density that depends on the maximal magnitude of the two variables: \begin{equation}
  f(x_1,x_2)={1 \over 2} \big(1-\Phi(\|\x\|_\infty)\big)
\end{equation}
where $\Phi(\cdot)$ is the cdf of standard normal distribution, and $\|\x\|_\infty= \max\{|x_1|,|x_2|\}.$  This bivariate density has a natural Bayesian interpretation: it can be regarded as the marginal density of $(X_1,X_2)$ if we put a uniform prior over the correlation $\rho$ (this type of density is referred to as the marginal predictive distribution in Bayesian literature). Moreover, one interesting feature of this density is that its isodensity contours consists of concentric squares.

Although we were not able to find the above result in the literature, we noticed that the bivariate density $f(x_1,x_2)={1 \over 2} \big(1-\Phi(\|\x\|_\infty)\big)$ is first obtained in a different manner by \citet{bryson1982}. In this paper, the authors consider constructing multivariate distributions through the Khintchine mixture method (\citet{khintchine1938}). The bivariate density $f(x_1, x_2)$ is listed as an example of their construction. But the link between this density and the uniform mixture over correlations is not addressed.

Through the Khintchine mixture approach, we show that the resulting mixed density is a function of $\|\x\|_\infty.$ Moreover, we show that for each $p$, this resulting density is the unique multivariate density that is a differentiable function of $\|\x\|_\infty$ and is marginally normal. It thus becomes interesting to investigate the connection between the Khintchine mixture and the uniform mixture over correlation matrices.

\section{The Uniform Correlation Mixture Integral}
Our first main result is the following theorem:
\begin{thm}\label{thm:gem}
  \begin{equation}
   f(x_1,x_2)=\int_{-1}^1{1 \over 2}f(x_1,x_2|\rho) d \rho ={1 \over 2} \big(1-\Phi(\|\x\|_\infty)\big).
  \end{equation}
\end{thm}
The proof can be found in Appendix A. Note that $f(x_1,x_2)$ is a proper bivariate density, and it is marginally standard normal:
\begin{equation}
\begin{split}
  & \int_\reals f(x_1,x_2) dx_2\\
  =& {1 \over 2}\int_{-|x_1|}^{|x_1|} \big(1-\Phi(|x_1|)\big) d x_2 + {1 \over 2}\int_{(-\infty,-|x_1|) \bigcup (|x_1|,\infty)} \big(1-\Phi(|x_2|)\big) d x_2\\
  =&\big(1-\Phi(|x_1|)\big)|x_1| + x_2\big(1-\Phi(x_2)\big) \bigg|_{|x_1|}^{\infty}+\int_{|x_1|}^\infty x_2 {1 \over \sqrt{2 \pi}} e^{-x_2^2/2} dx_2\\
  =&  \big(1-\Phi(|x_1|)\big)|x_1| - \big(1-\Phi(|x_1|)\big)|x_1| - {1 \over \sqrt{2 \pi}} e^{-x_2^2/2}\bigg|_{|x_1|}^\infty\\
  =&{1 \over \sqrt{2 \pi}} e^{-x_1^2/2}=\phi(x_1).
\end{split}
\end{equation}

The form of this bivariate density implies that its isodensity contours consists of concentric squares. Thus, an intuitive interpretation of this result is that if we ``average'' the isodensity contours of bivariate normal distributions which are concentrically elliptic, we shall get an isodensity contour of concentric squares. The plot of $f(x_1,x_2)$ is given in Figure 1.

\begin{figure}[hhhh]
\begin{center}
\includegraphics[width=4.5in]{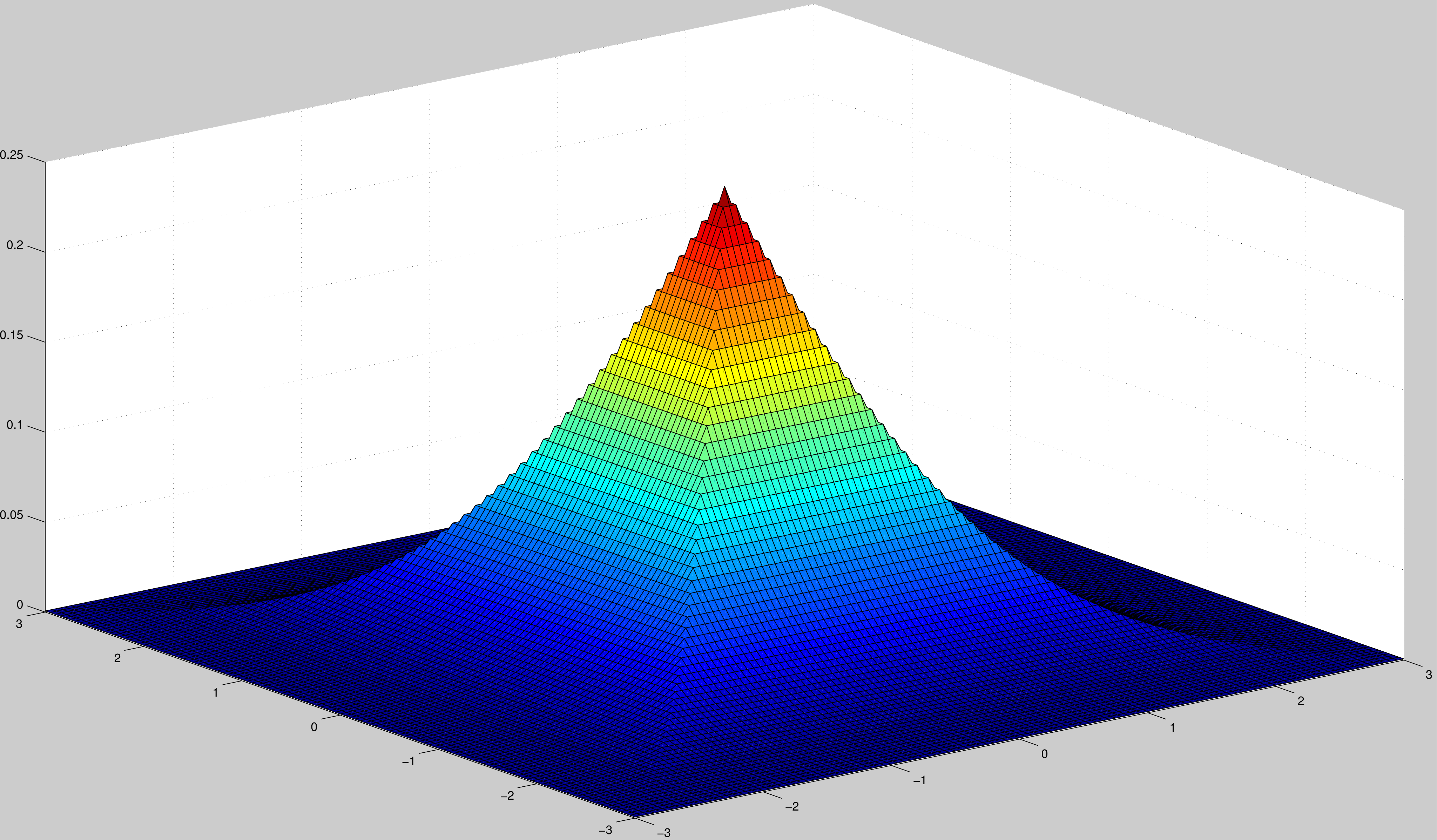}
\end{center}
\caption{The plot of the bivariate density function $f(x_1,x_2)={1 \over 2} \big(1-\Phi(\|\x\|_\infty)\big).$ Note that this bivariate density has a contour of squares and is marginally normal. }
\end{figure}

This result also indicates that if we have a uniform prior over the correlation $\rho,$ the resulting marginal density depends only on the maximal magnitude of the two variables. The application of this result to Bayesian inference needs to be further investigated. But from this marginal density of $x_1$ and $x_2$, the posterior distribution for $\rho$ can be derived in fully explicit form as the ratio of the bivariate normal density and this marginal density. Thus this result immediately leads to the Bayes factor of testing $\rho=0,$ in this very special situation. Moreover, the uniform prior has been used (see \citet{barnard2000} for theory and applications to shrinkage estimation) to model covariance matrices. We shall report these types of applications in future work.

\section{Connection to the Khintchine Mixture and Hypercubically Contoured Distributions that are Marginally Normal}
In \citet{bryson1982}, the authors developed a method of generating multivariate distributions through Khintchine's Theorem. Khintchine's theorem states that any univariate continuous random variable $X$ has a single mode if and only if it can be expressed as the product $X=Y U$ where $Y$ and $U$ are independent continuous variables and $U$ has a uniform distribution over $[0,1].$ This result and its extensions can be used to construct multivariate distributions that have specified marginal distributions. As an example of such a distribution with standard normal marginal distributions, the authors consider the construction with mutually independent $U_i \sim Uniform[-1,1], i=1,2$ and $Y \sim \chi_3$ (so that $Y>0$ and $Y^2 \sim \chi_3^2$). The random variables $X_1$ and $X_2$ are then generated as $X_1 =YU_1$ and $X_2=YU_2.$ The authors show that with this construction from $Y$ and $U_i$'s, the density of $(X_1,X_2)$ is exactly $f(x_1,x_2)= {1 \over 2} \big(1-\Phi(\|\x\|_\infty)\big).$

This density can also be generalized to higher dimensions through Khintchine's method. In fact, for any $p$, one generates $U_i \sim Uniform[-1,1], i=1,\ldots,p$ and $Y \sim \chi_3,$ and considers $X_i=Y U_i, i=1,\ldots,p.$ Then each $X_i$ is standard normally distributed. By using a $p+1$ dimensional transformation with $X_i=Y U_i$ and $W=Y$ and then integrating out $W$, we derive the joint density of $\X=(X_1,\ldots,X_p)$ as
\begin{equation}
  f_p(x_1,\ldots,x_p)={1 \over 2^{p-1} \sqrt{2 \pi}} \int_{\|\x\|_\infty}^\infty y^{2-p} e^{-y^2/2} dy.
\end{equation}
Note that $f_p(0,\ldots,0)=\infty$ for $p \ge 3.$ Nevertheless, $f_p$ is a proper $p$-dimensional density for every $p$ that has standard normal marginal distributions. Since $f_p$ is a function of $\|\x\|_\infty,$ the isodensity contour of $f_p$ consists of concentric hypercubes, which generalizes $f_2(x_1,x_2)={1 \over 2} \big(1-\Phi(\|\x\|_\infty)\big).$ We further show below that $f_p$ is the only density that possesses this property of being hypercubically contoured and having marginally normal distributions.

\begin{prop}\label{prop:unique}
  Consider a $p$-dimensional density that is a function of $\|\x\|_\infty,$ i.e., $g_p(x_1,\ldots,x_p)=h_p(\|\x\|_\infty)$ for some differentiable function $h_p: \reals^+ \rightarrow \reals^+.$ If $g_p(x_1,\ldots,x_p)$ has standard normal marginal densities, then the unique expression of $g_p(x_1,\ldots,x_p)$ is
  \begin{equation}\label{eq.unique}
    g_p(x_1,\ldots,x_p) = {1 \over 2^{p-1} \sqrt{2 \pi}} \int_{\|\x\|_\infty}^\infty y^{2-p} e^{-y^2/2} dy.
  \end{equation}
\end{prop}
The proof can be found in Appendix B.

\section{Discussion: Equivalence between the Uniform Correlation Mixture and the Khintchine Mixture in Higher Dimensions?}
In this paper, we have shown the equivalence of three bivariate densities: the uniform correlation mixture of bivariate normal densities, the unique square-contoured bivariate density with normal marginals, and the joint density of the Khintchine mixture of $\chi_3$ densities. We have also shown the equivalence of hypercubically-contoured densities and the Khintchine mixture of $\chi_3$ densities in higher dimensions. It thus becomes interesting to investigate whether the uniform correlation mixture of bivariate normal densities is equivalent to them in higher dimensions. Intuitively, this equivalence should carry over, as we would just be ``averaging'' the ellipsoid-shaped contours of multivariate normal densities instead of elliptical ones. However, directly integrating the multivariate normal density over the uniform measure over positive definite matrices is not a transparent task. Furthermore, if this relationship holds for normal distributions in higher dimensions, one may be curious whether this equivalence holds also for other distributions. The fact $f_p(0,\ldots,0)<\infty$ only for $p \le 2$ is also interesting. Due to the relationship between the normal distribution and random walks, it would be interesting to see if the finiteness of $f_p(0,\ldots,0)$ connects with the recurrence or transience of random walks.

Due to the natural Bayesian interpretation of this mixture, we will also investigate its Bayesian applications. For example, the Bayesian inference about $\rho$ discussed in this paper relies on the fact that the means and variances are known, and the correlation $\rho$ is the only parameter of interest. However, the important case with unknown means and variances (as studied in \citet{Geisser1964}) needs to be further studied, and it is not clear whether a fully explicit form of the posterior distribution exists in this case. We will study these problems in future work.

\appendix
\section{Proof of Theorem~2.1}
The derivation of $f(x_1,x_2)$ consists of the following steps:
\begin{enumerate}
  \item To show $f(0,0)={1 \over 4}={1 \over 2}\big(1-\Phi(0)\big).$
  \item To show that if $x_1$ and $x_2$ are not both $0$'s, then the function $g(\rho)= {x_1^2+x_2^2-2\rho x_1x_2 \over 2(1-\rho^2)}$ is monotone decreasing for $\rho \in (-1,a]$ and is monotone increasing for $\rho \in [a,1),$ where $a=a(x_1,x_2)=sgn(x_1x_2){|x_1|\wedge |x_2| \over |x_1| \vee |x_2|}.$
  \item To treat $I_1=\int_{-1}^a{1 \over 2}f(x_1,x_2|\rho) d \rho$ and $I_2=\int_{a}^1 {1 \over 2}f(x_1,x_2|\rho) d \rho$ separately, and to show that
  $$\int_{-1}^1 {1 \over 4 \pi \sqrt{1-\rho^2}} \exp \left\{-{x_1^2+x_2^2-2\rho x_1x_2 \over 2(1-\rho^2)}\right\} d \rho= I_1+I_2={1 \over 2} \big(1-\Phi(\|\x\|_\infty)\big).$$
\end{enumerate}

\noindent \textbf{Step 1}. For $x_1=x_2=0,$
\begin{equation}
f(0,0)= \int_{-1}^1  {1 \over 4 \pi \sqrt{1-\rho^2}} d \rho = {1 \over 4 \pi}\arcsin \rho\bigg|_{-1}^1={1 \over 4}={1 \over 2}\big(1-\Phi(0)\big).
\end{equation}

\noindent \textbf{Step 2}. If $x_1$ and $x_2$ are not both $0$'s, then consider the derivative of the exponent in $f(x_1,x_2|\rho)$: $$g(\rho)= {x_1^2+x_2^2-2\rho x_1x_2 \over 2(1-\rho^2)}.$$

After some algebra, it can be shown that
\begin{equation}
  {d \over d \rho} g(\rho) = -{1 \over (1-\rho^2)^2}(\rho x_1-x_2)(\rho x_2 - x_1).
\end{equation}
and
\begin{equation}
  {d^2 \over d^2 \rho} g(\rho) = {1 \over (1-\rho^2)^3}\big\{(1+3\rho^2)x_1^2-2\rho(\rho^2+3)x_1x_2+(1+3\rho^2)x_2^2\big\} \ge 0.
\end{equation}

Therefore, the minimum of $g(\rho)$ is attained at $a=a(x_1,x_2)=sgn(x_1x_2){|x_1|\wedge |x_2| \over |x_1| \vee |x_2|}.$ For example, if $x_2> -x_1 \ge 0,$ then $a={x_1 \over x_2}.$ We should also note that the minimum value of $g(\rho)$ is
\begin{equation}
  g(a)={|x_1|^2 \vee |x_2|^2 \over 2}.
\end{equation}

\noindent \textbf{Step 3}. Without loss of generality, we consider the case $x_1>x_2 \ge 0.$ We split the integral into two pieces:
\begin{equation}
  \int_{-1}^1{1 \over 2}f(x_1,x_2|\rho) d \rho = \int_{-1}^a {1 \over 2}f(x_1,x_2|\rho) d \rho + \int_{a}^1{1 \over 2}f(x_1,x_2|\rho) d \rho=I_1+I_2
\end{equation}
We start with
\begin{equation}
  I_2 = \int_a^1 {1 \over 4 \pi \sqrt{1-\rho^2}} \exp \left\{-{x_1^2+x_2^2-2\rho x_1x_2 \over 2(1-\rho^2)}\right\} d \rho.
\end{equation}
For $\rho \in [a,1),$ $g(\rho)$ is monotone increasing in $\rho$. Therefore, we consider the transformation $y=g(\rho).$ Solving the quadratic equation for $\rho$ yields that
\begin{equation}\label{eq.rho}
  \rho =g^{-1}(y)= {x_1x_2+\sqrt{(2y-x_1^2)(2y-x_2^2)} \over 2 y}.
\end{equation}
Denote $\Delta(y)=\sqrt{(2y-x_1^2)(2y-x_2^2)}.$ Note that
\begin{equation}
{d \over  d y} \Delta(y) = {4y-(x_1^2+x_2^2) \over \Delta(y)}.
\end{equation}
Thus, by differentiating \eqref{eq.rho} in $y$, we obtain
\begin{equation}
  d \rho = {(x_1^2+x_2^2)y-x_1^2x_2^2-x_1x_2 \Delta(y) \over 2 y^2 \Delta(y)} dy.
\end{equation}
Moreover, by \eqref{eq.rho},
\begin{equation}
  \sqrt{1-\rho^2} = {1 \over y}\sqrt{{1 \over 2}\bigg((x_1^2+x_2^2)y-x_1^2x_2^2-x_1x_2 \Delta(y) \bigg)}.
\end{equation}
By applying the last two equations to $I_2,$ we have
\begin{equation}
  I_2= \int_{x_1^2 \over 2}^\infty {\sqrt{(x_1^2+x_2^2)y-x_1^2x_2^2-x_1x_2 \Delta(y)} \over 4 \sqrt{2} \pi y \Delta(y)} e^{-y} dy.
\end{equation}

\noindent A similar argument for $I_1$ yields that
\begin{equation}
  I_1= \int_{x_1^2 \over 2}^\infty {\sqrt{(x_1^2+x_2^2)y-x_1^2x_2^2+x_1x_2 \Delta(y) } \over 4 \sqrt{2} \pi y \Delta(y)} e^{-y} dy.
\end{equation}

\noindent We show next that
\begin{equation}
\begin{split}
  & {\sqrt{(x_1^2+x_2^2)y-x_1^2x_2^2+x_1x_2 \Delta(y)} \over \sqrt{2(2y-x_2^2)}}+{\sqrt{(x_1^2+x_2^2)y-x_1^2x_2^2-x_1x_2 \Delta(y)} \over \sqrt{2(2y-x_2^2)}}\\
=& x_1
\end{split}
\end{equation}
\noindent This is seen by noting that
\begin{equation}
\begin{split}
  & \left({\sqrt{(x_1^2+x_2^2)y-x_1^2x_2^2+x_1x_2 \Delta(y)} \over \sqrt{2(2y-x_2^2)}}+{\sqrt{(x_1^2+x_2^2)y-x_1^2x_2^2-x_1x_2 \Delta(y)} \over \sqrt{2(2y-x_2^2)}}\right)^2\\
=& {4x_1^2y-2x_1^2x_2^2 \over 2(2y-x_2^2)}\\
=& x_1^2.
\end{split}
\end{equation}

\noindent Now the integral is simplified as
\begin{equation}
   \int_{-1}^1{1 \over 2}f(x_1,x_2|\rho) d \rho = I_1+I_2 =\int_{x_1^2 \over 2}^\infty {x_1 \over 4 \pi y \sqrt{2y-x_1^2}} e^{-y} dy = \int_{1 \over 2}^\infty {1 \over 4 \pi z \sqrt{2z-1}} e^{-x_1^2 z} dz
\end{equation}
where $z=y/x_1^2.$

\noindent By classical Laplace transform results in \citet{erdelyi1954}, the last integral is found to be
\begin{equation}
  \int_{1 \over 2}^\infty {1 \over 4 \pi z \sqrt{2z-1}} e^{-x_1^2 z} dz = {1 \over 2} \big(1-\Phi(x_1)\big).
\end{equation}

\section{Proof of Proposition~3.1}
We shall prove the proposition by induction. For $p=1,$ $\|\x\|_\infty=|x_1|,$ and the ``marginal'' density is the density in $x_1.$  Thus $g_1(x_1)=\phi(x_1)$ is the unique density that is ``marginally'' normal and is a function of $\|\x\|_\infty=|x_1|.$  This is the only case where the marginal normality is needed.

Suppose \eqref{eq.unique} is true for $1,\ldots p.$ For $p+1$, write $\tx=(x_1,\ldots,x_{p}).$ Suppose $g_{p+1}(\x)=h_{p+1}(\|\x\|_\infty)$ for some differentiable function $h_{p+1}: \reals^+ \rightarrow \reals^+$ and has standard normal marginal distributions.

Note that if we integrate out $x_{p+1},$ i.e.,
\begin{equation}
  \int_{-\infty}^\infty g_{p+1}(\x)d x_{p+1} = 2 \|\tx\|_\infty h_{p+1} (\|\tx\|_\infty) + 2 \int_{\|\tx\|_\infty}^\infty h_{p+1}(x_{p+1}) d x_{p+1},
\end{equation}
then this resulting density on the right hand side is a $p$-dimensional joint density that depends only on $\|\tx\|_\infty.$ By the uniqueness from the induction hypothesis,
\begin{equation}
  2 \|\tx\|_\infty h_{p+1} (\|\tx\|_\infty) + 2 \int_{\|\tx\|_\infty}^\infty h_{p+1}(x_{p+1}) d x_{p+1} = {1 \over 2^{p-1} \sqrt{2 \pi}} \int_{\|\tx\|_\infty}^\infty z^{2-p} e^{-z^2/2} dz.
\end{equation}

Now consider the equation in a positive $t$ that
\begin{equation}
  2th_{p+1}(t)+2\int_t^\infty h_{p+1}(x_{p+1}) d x_{p+1} = {1 \over 2^{p-1} \sqrt{2 \pi}} \int_t^\infty z^{2-p} e^{-z^2/2} dz.
\end{equation}
Differentiating in $t$ on both sides, we get
\begin{equation}
  h'_{p+1}(t)=-{1 \over 2^{p} \sqrt{2\pi}} t^{1-p} e^{-t^2/2}.
\end{equation}
Thus,
\begin{equation}
  h_{p+1}(y)= \int_{y}^\infty {1 \over 2^{(p+1)-1} \sqrt{2\pi}} t^{2-(p+1)} e^{-t^2/2} d t +c
\end{equation}
for some constant $c$. Since $\lim_{y\rightarrow \infty} h_{p+1}(y)=0,$ $c=0.$ Thus the induction proof is complete.

%% HERE WE DECLARE THE BIBLIOGRAPHYSTYLE TO USE AND THE BIBLIOGRAPHY DATABASE
\bibliographystyle{asa}
\bibliography{UCM_BiGaussian_1108}

\end{document}